\documentclass[a4]{amsart}
\usepackage{amssymb,amsmath,amsthm}
\usepackage{graphicx}
\usepackage{subfigure}
\usepackage{epsfig}
\usepackage[monochrome]{color}
\begin{document}
\title[Minimum polyhedron with $n$ vertices]{Minimum polyhedron with $n$ vertices}


\newtheorem{theorem}{Theorem}
\newtheorem{lemma}{Lemma}

\newtheorem{proposition}{Proposition}

\newtheorem{corollary}{Corollary}
\newtheorem{rem}[theorem]{Remark}
\newtheorem{ex}{Example}

\newcommand {\C}  {{\mathbb C}}
\newcommand {\N}  {{\mathbb N}}
\newcommand {\Z}  {{\mathbb Z}}
\newcommand {\R}  {{\mathbb R}}
\newcommand {\Q}  {{\mathbb Q}}
\renewcommand {\u} {{\mathbf{u}}}
\renewcommand {\v} {{\mathbf{v}}}
\newcommand {\w} {{\mathbf{w}}}
\newcommand {\x} {{\mathbf{x}}}
\newcommand {\y} {{\mathbf{y}}}
\newcommand{\Red}{\textcolor{red}}
\newcommand{\Blue}{\textcolor{blue}}

\subjclass{Primary 52B60, 52B55}
\keywords{}

\author[S. Akiyama]{Shigeki Akiyama}
\address{Institute of Mathematics\\
University of Tsukuba\\
1-1-1 Tennodai, Tsukuba, Ibaraki, 305-8571 JAPAN
} 
\email{akiyama@math.tsukuba.ac.jp}
\urladdr{http://math.tsukuba.ac.jp/$\sim$akiyama/}




\begin{abstract}
We study a polyhedron with $n$ vertices of fixed volume having minimum surface area.
Completing the proof of \Blue{Fejes T\'oth},
we show that all faces of a
minimum polyhedron are triangles, and further prove that a minimum polyhedron does 
not allow deformation of a single vertex.
We also present 
possible minimum shapes for $n\le 12$, some of them are quite unexpected, in particular $n=8$. 
\end{abstract}

\maketitle

\section{Introduction}
Let $X$ be a closed set in $\R^d$. 
Denote by $V_d(X)$  the $d$-dimensional Lebesgue measure of $X$ and
by $A_d(X)$ the $d-1$ dimensional Lebesgue measure of $\partial(X)$. 
For a non-empty set $A$ in $\R^d$, we denote by $\Delta(A)$ the convex hull
of $A$. 
A convex body in $\R^d$ is a compact convex set with a non-empty interior.
For a convex body $X$,  we recall the isoperimetric inequality:
\begin{equation}
\label{Iso}
\frac {A_d(X)}{V_d(X)^{(d-1)/d}}
\ge 
\frac {A_d(B^d)}{V_d(B^d)^{(d-1)/d}}
\end{equation}
where $B^d$ is the $d$-dimensional unit ball (c.f. \cite{Schneider:93}). 
The equality is attained only when $X$
is a $d$-dimensional ball. 
Note that if $X$ is a planar convex set, then in the plain language, 
$V_2(X)$ is the area and $A_2(X)$ is the perimeter of $X$. 

Let \Red{$d=3$ and} $n\ge 4$. We are interested in minimizing $A_3(X)/V_3(X)^{2/3}$ among
all polyhedra $X$ with $n$ vertices. Clearly we may assume
that $X$ is convex.  Denote by $\Delta_n=\Delta(p_1,\dots, p_n)$ 
the convex hull of $n$ points $p_1,p_2,\dots, p_n \in \R^3$. 
We say $\Delta_n$ is non-degenerate if $\Blue{V_3}(\Delta_n)>0$. 
Therefore our problem is to minimize $A_3(\Delta_n)/V_3(\Delta_n)^{2/3}$ 
among all non-degenerate convex hull $\Delta_n$'s of $n$ points in $\R^3$.
We are of course interested in the shape $\Delta_n$ which attains its minimum as well.
Clearly $A_3(\Delta_n)/V_3(\Delta_n)^{2/3}$ is invariant under similitudes, 
our problem is \Blue{equivalent to find} the minimum $A_3(\Delta_n)$ under $V_3(\Delta_n)=1$.
Thus our problem is a discrete variant of
the isoperimetric inequality (\ref{Iso}), i.e., a discrete `minimum surface'.
For a similar minimization problem with a given number of faces, we \Red{find}
several 
references. Lindel\"of \cite{Lindelof} and Minkowski \cite{Minkowski}
proved  in different intriguing ways that the minimum polyhedron 
must be circumscribed about a sphere, 
and \Red{Fejes T\'oth} \cite{Toth:48} proved that the minimum 
is attained when the number of faces are 4,6 and 12 by the regular tetrahedron, cube and dodecahedron, 
respectively. Note that minimization with a given number of vertices 
is a totally different and more difficult problem; e.g., 
the cube is not a solution for $n=8$ (see \Blue{Theorems} \ref{Triangle} and \ref{Minimum}).

\Red{Fejes T\'oth claimed in \cite{Toth:48_0} and \cite[Chapter 5, \S 7]{Toth} 
that every face of the minimum $n$-hedron must be a triangle.
However his proof contains a gap due to the fact that the corresponding 
equi-area piecewise differentiable surface 
has singular points (see the description around Example \ref{Singular} for details).
We shall classify such singularities (Lemma \ref{FiniteException}) and 
complete the proof along his idea in Theorem \ref{Triangle}.
All the same, digesting his idea, we can 
further prove that the minimum $n$-hedron does not allow
a deformation of a single vertex, by showing that the equi-area
body is strictly convex.
Finally we give a list of possible shapes of minimum $n$-hedron for
$n\le 12$ by extensive random numerical search together with heavy
algebraic computation using Gr\"obner basis.
Conjectural shapes for $n=8$ and $n=11$ may be beyond our imagination.
}

\Red{Added in revision: One of the referees of this paper pointed out that
B\"{o}r\"{o}czky and B\"{o}r\"{o}czky~Jr \cite{Boroczky2_Iso:97} gave a proof of Theorem 
\ref{Triangle} and its generalization.
Their proof also rely on the same idea of Fejes T\'oth but went in a different way. 
}

\section{\Red{Every face is a triangle}}
\Red{
In this section, we prepare 
basic properties of this minimization problem in Propositions \ref{Exist} and \ref{Decrease}.
Then we point out a gap in the proof of Fejes T\'oth \cite{Toth:48_0, Toth} which asserts that 
every face of the minimum $n$-hedron is triangular. Then we
complete the proof after the classification of singular points of the 
equi-area surface.}

\begin{lemma}
\label{OProj}
Let $Y$ be a planar polygon in $\R^3$ and $g:\R^3\rightarrow \R^2$ be an orthogonal 
projection to some plane (for e.g., the one along $z$-axis to $xy$-plane). 
Then we have $V_2(g(Y))\le V_2(Y)$ and $A_2(g(Y))\le A_2(Y)$.
\end{lemma}

\proof
This is clear from the property $\| g(x)-g(y)\|\le \|x-y\|$
for any $x,y$. 
\qed
\medskip

\begin{proposition}
\label{Exist}
For a fixed integer $n\ge 4$,  the minimum of  
$A_3(\Delta_n)/V_3(\Delta_n)^{2/3}$ exists 
where $\Delta_n$ varies among non-degenerate
convex hulls of $n$ points in $\R^3$.
\end{proposition}

\proof
Let $R$ be the diameter of $\Delta(p_1,\dots, p_n)$ attained by $\| p_1-p_2\|=R$. Let $S$ be the plane passing 
through $p_1$ which is 
orthogonal to the segment $[p_1,p_2]$ and $g$ be the orthogonal projection to $S$.
Then $g(\Delta_n)$ is a convex polygon in $S$ with 
vertices $q_1,\dots,q_{\ell}$ with $\ell\le n-1$,  arranged in the clockwise order
with respect to the centroid of $g(\Delta_n)$.
Choose $q'_1,\dots, q'_{\ell}$ 
in $\Delta_n$ such that $g(q'_i)=q_i$ for $i=1,\dots,\ell$. 
\Red{If the segment $[p_1,p_2]$ lies within $\partial(\Delta(p_1,\dots, p_n))$,
we choose $q'_1=(p_1+p_2)/2$.}
Since $\Delta_n$ is contained in $g(\Delta_n) \times [0,R]$ we have
\begin{equation}
\label{V}
V_3(\Delta_n)\le V_2(g(\Delta_n)) R.
\end{equation}
We claim that 
\begin{equation}
\label{A}
A_3(\Delta_n)\ge \frac 12 A_2(g(\Delta_n)) R.
\end{equation}
\Red{
Since
$X\supset Y$ implies $A_3(X)\ge A_3(Y)$ for convex bodies $X, Y$,
considering the surface area of the convex hull $Y=\Delta(
p_1,p_2, q'_1,q'_2,\dots, q'_{\ell})$, it is enough to prove} that $$
V_2(p_1,q'_i,q'_{i+1})+
V_2(p_2,q'_i,q'_{i+1}) \ge \frac 12 \| q_i-q_{i+1}\| R,$$
where $V_2(x,y,z):=V_2(\Delta(x,y,z))$, the area of the triangle of vertices $x,y,z$. 
The index $i$ of $q_i$ is considered modulo $\ell$. 
Take a plane $P$ 
containing $p_1$ and $p_2$ parallel to the segment $[q'_i,q'_{i+1}]$ and use
the orthogonal projection $g_2$ to $P$. \Red{Noting} that the directions of the two projections 
$g$ and $g_2$ are orthogonal, we have $g(q'_i)-g(q'_{i+1})=g(g_2(q'_i))-g(g_2(q'_{i+1}))$.
By Lemma \ref{OProj}, \Red{we see
$$
V_2(p_1,q'_i,q'_{i+1})+V_2(p_2,q'_i,q'_{i+1}) \ge \Blue{
V_2(p_1,g_2(q'_i),g_2(q'_{i+1}))+V_2(p_2,g_2(q'_i),g_2(q'_{i+1}))}$$
which is not less\footnote{\Red{This holds even when $[g_2(q'_i),p_1]$ and $[g_2(q'_{i+1}),p_2]$ intersect.}}}
than 
$\|g(q'_i)-g(q'_{i+1})\| R/2$. This shows the claim. 

Using (\ref{V}) and the 
isoperimetric inequality (\ref{Iso}) for $d=2$, that is,  $A_2(g(\Delta_n))^2\ge 4\pi V_2(g(\Delta_n))$ ,
we see $A_2(g(\Delta_n))\ge 2 \sqrt{\frac {\pi V_3(\Delta_n)}R}$. 
Let us fix $V_3(\Delta_n)=1$. 
Then we have $A_3(\Delta_n)\ge \sqrt{\pi R}$ from (\ref{A}). This shows that
$A_3(\Delta_n)\rightarrow \infty$ as $R\rightarrow \infty$ under the assumption
$V_3(\Delta_n)=1$. Since we are interested in minimizing $A_3(\Delta_n)/V_3(\Delta_n)^{2/3}$, 
we may assume that $R$ is bounded by some constant $K$. 
This shows that
parameters $p_1,\dots, p_n$ are in a closed 
ball of radius $K$ with the prescribed property $V_3(\Delta_n)=1$. Therefore the set of
parameters are in a compact set in $\R^3$ 
and we find the minimum of $A_3(\Delta_n)$ as desired.
\qed
\medskip

Therefore we define $\alpha_n=\min_{\Delta_n} A_3(\Delta_n)/V_3(\Delta_n)^{2/3}$ where $\Delta_n$ runs over all
non-degenerate convex hulls of $n$ points. A {\it minimum $n$-hedron} is 
the shape $\Delta_n$ which attains $\alpha_n$. It may not be unique but we expect that 
it is unique up to similitudes in $\R^3$.

\begin{proposition}
\label{Decrease}
We have 
$\alpha_n>\alpha_{n+1}$ for $n\ge 4$ and $\lim_{n\rightarrow \infty}\alpha_{n}=(36\pi)^{1/3}\approx 4.83598$.
\end{proposition}

\proof
Choose $\Delta_n$ which attains $\alpha_n$ and 
its face $T\subset \partial(\Delta_n)$. 
We take a point $p_{n+1}$ 
on a outward normal emanating from an inner point $p$ of $T$ whose 
distance from $T$ is $\varepsilon>0$, which is small enough that $\Delta_{n+1}$ is 
the union of $\Delta_n$ and the pyramid of base $T$ and the apex $p_{n+1}$.
Denote by $e_i$ the edge of $T$ and $r_i$ be the height of $p$ from the edge $e_i$
for $i=1,2,\dots, t$. Note that $r_i>0$.
Then we see
$$V_3(\Delta_{n+1})=V_3(\Delta_n)+\frac 13 \varepsilon V_2(T)$$
and
$$A_3(\Delta_{n+1})=A_3(\Delta_n)-V_2(T)+\frac 12\sum_{i=1}^t e_i \Red{\sqrt{\varepsilon^2+r_i^2}.}$$
Using $V_2(T)=\frac 12 \sum_{i=1}^t e_i r_i$, we have
\begin{eqnarray*}
\frac {A_3(\Delta_{n+1})}{V_3(\Delta_{n+1})^{2/3}}&=&
\frac {A_3(\Delta_{n})}{V_3(\Delta_{n})^{2/3}} 
\frac {1+\frac {1}{2 A_3(\Delta_{n})} \sum_{i=1}^t e_i r_i \left(\sqrt{1+(\frac{\varepsilon}{r_i})^2} -1\right)}
{\left(1+\frac {\varepsilon V_2(T)}{3V_3(\Delta_n)} \right)^{2/3}}\\
&=&
\frac {A_3(\Delta_{n})}{V_3(\Delta_{n})^{2/3}} \frac {1+C_1\varepsilon^2+O(\varepsilon^3)}
{1+C_2\varepsilon+O(\varepsilon^2)}
\end{eqnarray*}
with $C_1=\frac 1{4A_3(\Delta_n)} \sum_{i=1}^t \frac {e_i}{r_i}$ and $C_2=\frac {2V_2(T)}{9V_3(\Delta_n)}$. 
Taking \Blue{small $\varepsilon>0$,} we have
$$
\alpha_n=\frac {A_3(\Delta_{n})}{V_3(\Delta_{n})^{2/3}}>
\frac {A_3(\Delta_{n+1})}{V_3(\Delta_{n+1})^{2/3}} \ge \alpha_{n+1}.
$$
By isoperimetric inequality (\ref{Iso}) for $d=3$, we have
$$
A_3(\Delta_n)/V_3(\Delta_n)^{2/3}\ge (36\pi)^{1/3}
$$
and the minimum is sufficiently approximated by points on the sphere, provided
$n$ is large.
\qed
\medskip

\begin{theorem}
\label{Triangle}
Every face of a minimum $n$-hedron is a triangle. 
\end{theorem}

The statement is intuitively quite natural, because we want a round shape and 
bending non-triangular faces by pulling outward their diagonals 
does not increase the number of vertices. 
We shall prove Theorem \ref{Triangle} after Lemma \ref{FiniteException}.  
Here we quote a paragraph on Theorem \ref{Triangle} 
in page 58 of \Red{Fejes T\'oth} \cite{Toth:48_0} (see also \cite{Toth:39, Toth}). 
\medskip

{\it
Greifen wir um dies einzusehen eine beliebige Ecke E des als extremal vorausgesetzten Polyeders heraus und bewegen es so, dass erstens der Inhalt, zweitens die \Red{Oberfl\"ache} der kleinsten konvexen H\"ulle $\mathbb{H}$ von E und der \"Ubrigen Ecken des Polyeders konstant bleiben. Im ersten Fall durchl\"auft E den Rand eines konvexen Polyeders $\mathbb{P}$, im zweiten Fall dagegen den Rand eines singularit\"atenfreien Eik\"orpers $\mathbb{E}$, der im Falle \Red{eines} Extremalen Polyeders offenkundig keinen Punkt ausserhalb $\mathbb{P}$ \Red{haben} kann. W\"ahre nun E die Ecke einer mehr als dreiseitigen Fl\"ache des urspr\"unglichen Polyeders, so liege E -wie eine einfache \"Uberlegung zeigt- auf einer Kante von $\mathbb{P}$. Mithin k\"onnte $\mathbb{P}$ nicht die singularit\"atenfreie Fl\"ache $\mathbb{E}$ enthalten.
}
\medskip

{\it \Red{(English Translation)}
In order to see this,  let us take an arbitrary vertex E of the polyhedron, which is supposed to be 
extremal, 
we can move it keeping firstly the volume, and secondly, the surface area of 
the smallest convex hull $\mathbb{H}$ of E and the remaining vertices of the polyhedron. 
In the first case, E goes through the boundary of a convex polyhedron $\mathbb{P}$,
and in the second case the boundary of a singularity-free body $\mathbb{E}$, which, 
in the case of an extreme polyhedron, is obviously not a point outside $\mathbb{P}$. 
If E is the vertex of a more than three-sided face of the polyhedron, then, 
by a simple discussion, it is on an edge of $\mathbb{P}$. 
However $\mathbb{P}$ can not contain the singularity-free surface $\mathbb{E}$.
}
\medskip

Let us try to understand this description \Red{and then show where is the gap}.
A point $\x$ in a convex set $X$ is {\it visible} from a point $\y\in \R^3\setminus X$,
if the segment $[\x,\y]$ intersects $X$ only at $\x$. 
A subset $V$ of $X$ is visible from $\y$ if each element of $V$ is visible from $\y$.
A {\it face plane} of a polyhedron 
$X$ is a hyperplane containing a \Red{codimension} one face of $X$. 
\Red{Let us fix $p_1,p_2,\dots p_{n-1}$.}
Define
\begin{equation}
\label{C}
C=\{ \v\in \R^3|\ V_3(\Delta(p_1,\dots,p_{n-1},\v))\le v\}
\end{equation}
and
\begin{equation}
\label{S}
S=\{ \v\in \R^3|\ A_3(\Delta(p_1,\dots,p_{n-1},\v))\le h\}
\end{equation}
with $v>V_3(\Delta(p_1,\dots,p_{n-1}))$ and $h>A_3(\Delta(p_1,\dots,p_{n-1}))$.
Clearly $V_3(\Delta(p_1,\dots,p_{n-1},\v))$ and $A_3(\Delta(p_1,\dots,p_{n-1},\v))$ are
continuous functions of $\v$.
The boundary $\partial(C)$ is a contour 
of the volume function of convex hull of $\v$ and visible faces from $\v$
of  $\Delta(p_1,\dots,p_{n-1})$. 
Visible faces change when and only when 
$\v$ passes a face plane of $\Delta(p_1,\dots,p_{n-1})$ and that
makes a visible face $F$ to a non-visible one, or the other way round. 
Note that 
this change happens only when a non-triangular face with a vertex $\v$ 
appears in $\Delta(p_1,\dots,p_{n-1},\v)$. 
This volume is an affine function 
on the coordinates of $\v$ determined by visible faces from $\v$. 
Therefore $C$ is the intersection of half-spaces defined by visible faces, i.e., a convex polyhedron. 
If there exists a non-triangular face, then $\v$ must be on the edge of $\partial(C)$.

The surface $\partial(S)$ is determined by visible edges 
from $\v$ which contribute the surface of the convex hull. 
Locally $\partial(S)$ is defined as a contour of sum of 
square roots of quadratic polynomials of its coordinates, 
which implies that $\partial(S)$ is piecewise smooth.
While $\v$ moves around, visible edges will switch to new ones
when the visible faces change. Note that this change happens when $\v$ is
on a face plane of $\Delta(p_1,\dots,p_{n-1})$. 
Assume that $\Delta_n=\Delta(p_1,\dots,p_{n-1},p_n)$ is a minimum $n$-hedron \Red{
and $p_n$ is a vertex of a non-triangular face. Then $p_n$ must be on the edge of $\partial(C)$.} 
Put
$v=V_3(\Delta_n)$ and $h=A_3(\Delta_n)$.
If $p_n$ is a totally differentiable point of $S$, then the surface $\partial(S)$ penetrates $\partial(C)$ 
and we must have a point of $\partial(S)$ outside $C$, which contradicts the minimality of $\Delta(p_1,\dots,p_n)$.
Therefore in this view, if $\partial(S)$ is totally differentiable everywhere, the proof is done (see 
\cite[Lemma 4.4]{Boroczky2_Iso:97}).

This idea is very insightful but does not work as it is. 
Here is a counter example \Red{for $n=5$} \Blue{such}
that $\partial(S)$ has a singular point. 

\begin{ex}
\label{Singular}
Let $p_1=(0,1,0),p_2=(0,0,0),p_3=(1,0,0),p_4=(0,0,1)$. The surface $S$ with
\Red{$h=4$ consists of $14$ algebraic surfaces and the black point 
$(0,2,0)$ is a common point of \Blue{four} surfaces} depicted in Figure \ref{Sing}.
\begin{figure}[htbp]
\begin{center}
\includegraphics{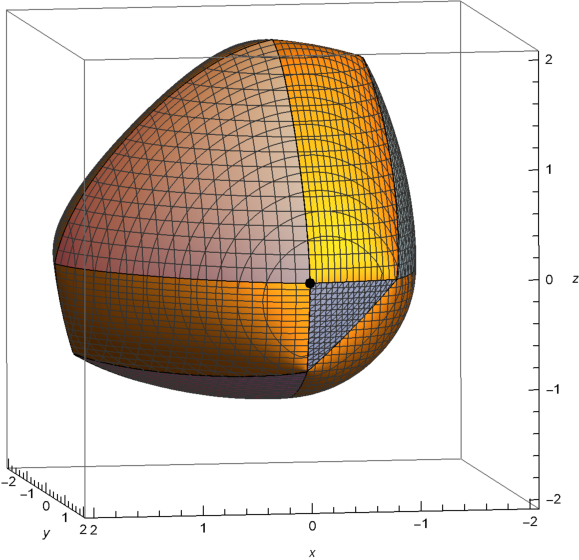}
\end{center}
\caption{The surface $\partial(S)$ with singular points}
\label{Sing}
\end{figure}
We claim that \Red{this point} $(0,2,0)$ is singular.  Indeed it is on the boundary of four surfaces
\begin{equation}
\label{1Q}
\sqrt{2 x^2+(y+z-1)^2}+\sqrt{(x+z-1)^2+2 y^2}+\sqrt{(x+y-1)^2+2 z^2}+3=8,
\end{equation}
$$
\sqrt{x^2+y^2}+\sqrt{(x+z-1)^2+2 y^2}+\sqrt{y^2+z^2}+1=8,
$$
$$
\sqrt{2 x^2+(y+z-1)^2}+\sqrt{x^2+z^2}+\sqrt{(x+z-1)^2+2 y^2}+\sqrt{y^2+z^2}+2=8,
$$
$$
\sqrt{x^2+y^2}+\sqrt{x^2+z^2}+\sqrt{(x+z-1)^2+2 y^2}+\sqrt{(x+y-1)^2+2 z^2}+2=8
$$
whose domains and visible edges \Blue{which contribute to $A_3(\Delta(p_1,p_2,p_3,p_4,\v))$} are 
$$\{(x,y,z)|\ x\ge 0, y\ge 0, z\ge 0\Red{, x + y + z \ge 1}\},\Red{\{\{p_4,p_1\},\{p_3,p_4\},\{p_1,p_3\}\}}$$ 
$$\{(x,y,z)|\ x\le 0, y \ge 0, z\le 0 , x + y + z \ge 1\},
\Red{\{\{p_2,p_4\},\{p_4,p_3\},\{p_3,p_2\}\}}$$
$$\{(x,y,z)|\ x\ge 0, y\ge 0, z\le 0,x+y+z\ge 1\}, \Red{\{\{p_4,p_1\},\{p_2,p_1\},\{p_3,p_4\},\{p_3,p_2\}\}}$$ 
$$\{(x,y,z)|\ x \le 0, y \ge 0, z\ge 0 , x + y + z \ge 1\}, 
\Red{\{\{p_2,p_4\},\{p_2,p_1\},\{p_3,p_4\},\{p_3,p_1\}\}}$$
\Red{respectively. For example, (\ref{1Q}) follows from
\begin{align*}
&V_2(p_4, p_1, \v) + V_2(p_3, p_4, \v) +  V_2(p_1, p_3,\v)\\
&+ V_2(p_1, p_2, p_3)+V_2(p_1, p_2, p_4)) +V_2(p_2,p_3,p_4)=4.
\end{align*}}
For the first two surfaces, 
outer normals at $(0,2,0)$ approaching from the corresponding
domains are $(1,5,1),(-1,10,-1)$, which are mutually 
inconsistent and tangent plane at $(0,2,0)$ can not be defined. 
For the \Red{remaining two surfaces}, the situation is worse that $(0,2,0)$ becomes a singular point by \Red{the} effect of
the term $\sqrt{x^2+z^2}$, whose partial derivatives on $x,z$ varies by the ratio $x:z$. 
\Red{We shall see in Lemma \ref{FiniteException}
that this type of singularity never vanishes regardless of the choice of $h$.}
\end{ex}

A {\it polyhedral cell} is a 
closed convex set with a non-empty interior whose boundaries consist of 
finite number of convex subsets of hyperplanes of \Red{codimension} 1. 
Hereafter we use a partition of $\R^3$ into polyhedral cells
by face planes $W_i$ of a convex hull $\Delta$. 
For $\v\in \R^3\setminus \Delta$, consider a plane $W$ separating $\Delta$ and $v$.
Then the union of visible faces from $\v$ is homeomorphically mapped to a figure of $W$ by a 
projection sending a point $y$ on the union to the point $y'\in W$ if $y,y',\v$ are collinear. 
We say
that the resulting figure is the {\it planar projection}. Planar projections are affine equivalent under the
change of separating planes. We prepare an important property of visibility.

\begin{lemma}
\label{Proj}
The planar projection of the union of visible faces $\{Q_i\}_{i=1}^k$ from $\v$ is convex. 
\end{lemma}

\proof
This follows immediately from the convexity of $\Delta$. 
\qed
\medskip

We first confirm that \Red{Fejes T\'oth}'s idea is almost valid, however, 
the surface $\partial(S)$ must have a singular point.

\begin{lemma}
\label{FiniteException}
Assume that $\Delta(p_1,\dots,p_{n-1})$ is non-degenerate and 
fix a positive constant $h>A_3(\Delta(p_1,p_2\dots,p_{n-1}))$. 
The surface $$\partial(S)=\{\v\in \R^3 |\ A_3(\Delta(p_1,\dots,p_{n-1},\v)) =h \}$$ is totally differentiable 
except at most $2e$ points where $e$ is the number of edges of $\Delta(p_1,\dots,p_{n-1})$.
The surface $\partial(S)$ is not totally differentiable at $\v\in \partial(S)$ if and only if the prolongation 
of an edge of $\Delta(p_1,\dots,p_{n-1})$ penetrates $\v$.
\end{lemma}

\begin{figure}[htbp]
\begin{center}
\includegraphics[height=4cm]{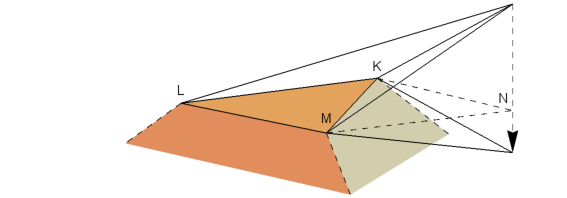}
\end{center}
\caption{Switching visible edges}
\label{Switch}
\end{figure}

\proof
We prove that $\partial(S)$ is totally differentiable at a switching point contained
in exactly one face plane of $\Delta(p_1,\dots,p_{n-1})$. 
This switching occurs at several 
contiguous edges forming a broken line which are the edges of the changing face. 
A crucial point is that the initial and the final vertex of this broken 
line does not change by the switching.
For example, consider a (planar) convex quadrangle $KLMN$ with $K=p_1,L=p_2,M=p_3$.
This could be divided into two triangles 
in two different ways, like $KLN$, $LMN$ or $KLM$, $KNM$. 
Edge switching occurs when $\v$ passes \Red{transversally} through $N$. 
First the area of triangles $KL\v$ and $LM\v$ contributes to $A_3$ 
 and later, triangle $KLM$ and $K\v M$ does. 
In this case, the related edges are $KL$, $LM$ at the beginning and switched to $KM$ (see Figure \ref{Switch}).

To see that $\partial(S)$ admits a tangent plane at the switching point $\v$, assume that 
the switching of edges happens on the plane $z=0$ to simplify the computation.
Let
$$(a_i,b_i,0),(a_{i+1},b_{i+1},0) \qquad (i=1,\dots,k-1)$$ 
be the end points
of the switching edges and $(x,y,z)\in \partial(S)$. Then surface $\partial(S)$ \Red{
on one side of the switching plane} 
is defined locally by an equation of the form $f(x,y,z)+m(x,y,z)=h$:
$$
f=\sum_{i=1}^{k-1} \frac{
\sqrt{ \left((b_i-b_{i+1})^2 + (a_i-a_{i+1})^2\right)z^2 +
((b_i-b_{i+1})x-(a_i-a_{i+1})y+a_ib_{i+1}-b_ia_{i+1})^2}}2
$$
where $m(x,y,z)$ is the contribution from non switching edges.
\Red{The equation and the value $k$ change when $v$ passes
the switching plane.}
At a switching point $(x_0,y_0,0)$, we have 
$$
\frac{\partial{(f+m)}}{\partial x}
=
\sum_{i=1}^{k-1} \frac{b_i-b_{i+1}}2 \frac
{(b_i-b_{i+1})x_0-(a_i-a_{i+1})y_0+a_ib_{i+1}-b_ia_{i+1}}
{\sqrt{((b_i-b_{i+1})x_0-(a_i-a_{i+1})y_0+a_ib_{i+1}-b_ia_{i+1})^2}}
+\frac{\partial{m}}{\partial x},
$$
$$
\frac{\partial{(f+m)}}{\partial y}
=
\sum_{i=1}^{k-1} \frac{a_{i+1}-a_{i}}2 \frac
{(b_i-b_{i+1})x_0-(a_i-a_{i+1})y_0+a_ib_{i+1}-b_ia_{i+1}}
{\sqrt{((b_i-b_{i+1})x_0-(a_i-a_{i+1})y_0+a_ib_{i+1}-b_ia_{i+1})^2}}
+\frac{\partial{m}}{\partial y}
$$
and
$$
\frac{\partial{(f+m)}}{\partial z}
=\frac{\partial{m}}{\partial z}.
$$
From Lemma \ref{Proj}, we see that $(a_i,b_i,0)\ (i=1,\dots,k)$ and $(x_0,y_0,0)$
\Red{form vertices of a planar} convex $(k+1)$-gon. Consequently
$$(b_i-b_{i+1})x_0-(a_i-a_{i+1})y_0+a_ib_{i+1}-b_ia_{i+1}$$ have the same sign
for all $i$ and the normal vector of $\partial(S)$ at $(x_0,y_0,0)$ 
is $(b_1-b_{k},a_k-a_1,0)\pm 2m'$, where $m'$ is 
the contribution from $m(x,y,z)$.
The vector $(b_1-b_{k},a_k-a_1,0)$ is 
orthogonal to the segment joining two end points of the switching broken line, i.e., 
the segment between $(a_1,b_1,0)$ and $(a_k,b_k,0)$. 
As these two end points are invariant under switching, 
even at the switching point $(x_0,y_0,0)$ 
the tangent plane is well-defined. 
Therefore the surface $\partial(S)$ is totally differentiable at any switching point
contained in exactly one face plane of $\Delta(p_1,\dots,p_{n-1})$. 

Let us study the possible singularities. 
The switching points lying on two or more face planes
are on the intersection lines of face planes. 
Assume that a point $\v$ is lying on face planes $W_j\ (j=1,\dots,\ell)$ with $\ell\ge 2$. 
The face plane 
$W_j$ induces switching of edges $e_1^{(j)},\dots, e_{k_j}^{(j)}\ (k_j\ge 2)$ 
to $e'^{(j)}$ or the other way round.\footnote{In the above proof,
the face plane is $z=0$ and
$e_i^{(j)}=[(a_i,b_i,0),(a_{i+1},b_{i+1},0)], e'^{(j)}=[(a_1,b_1,0),\Red{(a_{k},b_{k},0)}]$.}  
Of course $e'^{(j)}\ (j=1,\dots,\ell)$ are distinct. If the set of edges
$E_j=\{e_i^{(j)}|\ i=1,\dots, k_j\}$ are mutually disjoint for $j=1,\dots,\ell$, then $\partial(S)$
is totally differentiable at $\v$ by the same proof. 
The singularity happens
only when there exist $j_1\neq j_2$ \Red{for which} $E_{j_1}\cap E_{j_2}\neq \emptyset$. 
We claim that this is also sufficient.  In fact,
such an intersection must be a single edge and its prolongation must pass through the point $\v$. 
This means that around $\v$, there is a polyhedral cell 
$K$ that if $\u\in K$ then
there is a visible edge from $\u$ penetrating $\v$, that 
contributes the sum of the surface area of $\Delta(p_1,\dots,p_{n-1},\u)$. This contribution
is the square root of a positive \Red{semi-}definite 
quadratic form over three variables $x-a,y-b,z-c$ with $\v=(a,b,c)$,
which vanishes\footnote{In Example \ref{Singular}, this is the term $\sqrt{x^2+z^2}$
\Red{which vanishes on the line through the edge $\{p_1,p_2\}$.}}
\Red{when and only when 
$\u$ is on the line passing the visible edge penetrating $\v$}.
\Red{
Such a term corresponds exactly to two visible faces and gives a conic singularity. In most cases, such a term is unique and 
it produces a singularity at $\v$.
\Blue{For special cases,} 
there may be several terms of this type, 
each term vanishes on different 
prolongations of edges of $\Delta(p_1,\dots,p_{n-1})$
penetrating $\v$.
If there are more than one such terms, then all partial derivatives of the 
terms at $\v$ with respect to variables $x-a,y-b,z-c$ are zero.
Therefore the singularity at $\v$ becomes removable after the summation
only when this algebraic function becomes locally a constant.
However then the sum is constant everywhere, by algebraicity. This
cannot happen because every term diverges to positive infinity by 
taking limit in all directions except the vanishing line.
Therefore $\partial(S)$ cannot be totally differentiable at $\v$.}
This shows the claim and finishes the proof.
\qed
\medskip

{\it Proof of Theorem \ref{Triangle}.}
If there exists a non-triangular face and \Red{$p_n$ is a vertex on this face}, then 
\Red{$\v=p_n$} must be on an edge of $\partial(C)$. \Red{Here $C$ and $S$ are defined by (\ref{C}) and (\ref{S}) respectively}. 
\Red{By the discussion before Example \ref{Singular}, if $\v=p_n$ is non-singular then we get a contradiction.
Let $\v=p_n$ be} one of the singular points of $\partial(S)$ 
in Lemma \ref{FiniteException}. \Red{Then} 
there is an edge of  $\Delta(p_1,\dots,p_{n-1})$  penetrating $p_n$. 
However this implies that one of the \Red{vertices} of $\Delta(p_1,\dots,p_{n-1})$ is 
in the relative interior of an edge of $\Delta(p_1,\dots,p_{n-1},p_n)$. 
In this case the number of vertices of $\Delta(p_1,\dots,p_{n})$ is less than $n$.
Since $\Delta(p_1,\dots,p_{n-1},p_n)$ is a minimum
$n$-hedron, this does not happen by Proposition \ref{Decrease}, giving a contradiction \Red{for this case}.
\qed
\medskip

It is possible to give a geometric (but more technical)  alternative
proof of Theorem \ref{Triangle} without using the last characterization of the
singularity in Lemma \ref{FiniteException}. 
We give a rough sketch of it.
The singular point of $\v=p_n$ of $\partial(S)$
in the above proof is defined by piecewise smooth surfaces.  Take a polyhedral cell
$K_1$, $K_2$
defined by the face planes  $W_j\ (j=1,\dots \ell)$ with maximum and minimum number of
visible faces. It is clear that on $K_i\ (i=1,2)$ we see no visible edges passing $\v$,  and
therefore the tangent planes approaching from $K_i$ are well defined. This tangent plane
must coincide with the corresponding face planes of $\partial(C)$ (otherwise one can prolong a
tangent plane which penetrates $\partial(C)$ giving a smaller $A_3(\Delta)/V_3(\Delta)^{2/3}$ 
by a non convex $\Delta$). On the other hand approaching to $\v$ from 
other polyhedral cells surrounding $\v$,
the point $\v$ is singular. Partial derivatives of the singular terms appear in this intermediate terms
are the function on the ratio of $x-a:y-b:z-c$ with $\v=(a,b,c)$. Take a slice of $\partial(S)$ by a plane
which passes an inner point of $C$ close to $\v$ and intersects all $W_j$. 
This gives a piecewise smooth planar curve
that has two `almost' linear parts and other parts with positive curvature. Shifting 
the slice plane parallel and
closer to $p$, the shape converges to a single curve up to similitude, which encircles a convex planar 
region. On the other hand, since tangent planes exist within $K_i$, the parts of 
the curve in $K_i$ converge to line segments. Recalling $S\subset C$, 
this causes an inconsistency at their end points.
\medskip

\section{\Red{One vertex deformation is impossible}}
\Red{In this section, we further develop 
the idea of Fejes T\'oth and show that the minimum $n$-hedron does not allow
a deformation of a single vertex.}

Let $X$ be a convex set in $\R^d$. 
A function $F:X\rightarrow \R$ is convex if for any $\u,\v\in X$ and any $\lambda\in [0,1]$, we have
\begin{equation}\label{Conv}
F((1-\lambda)\u+\lambda \v)\le (1-\lambda)F(\u)+\lambda F(\v). 
\end{equation}
It is {\it strictly convex}
if for any $\u,\v\in X$ with $\u\neq \v$ and any $\lambda\in (0,1)$, 
$$
F((1-\lambda)\u+\lambda \v)<(1-\lambda)F(\u)+\lambda F(\v).
$$
Take a convex subset $Y\subset X$. 
If $F:X\rightarrow \R$ is convex and
the equality of (\ref{Conv}) with $\lambda\in (0,1)$ holds 
only when $\u,\v\in Y$, then we say $X$ is strictly convex except $Y$.

The next lemma gives a method to paste together convex functions defined in polyhedral cells to obtain a global convex function. Related general criteria are found in \cite{BauschkeYvesHung:16} using convex analysis. 

\begin{lemma}
\label{Piece}
Let $\R^d$ be partitioned into a finite number of polyhedral cells $\{D_i\}$ whose
interiors are disjoint.
Let $Z$ be the set of points of $\R^d$ that belong to more than two $D_i$.
Assume that $F_i$ is a convex function on $D_i$ so that $F_i(\v)=F_j(\v)$ holds for each 
$\v\in D_i\cap D_j$.
Then the function $F:\R^d \rightarrow \R$ is naturally defined by the values of $F_i$. 
We see that
$F$ is convex if and only if the following condition holds
\begin{itemize}
\item If $\v\in (D_i\cap D_j)\setminus Z$, 
$\v-\omega\in D_i$ and $\v+\omega\in D_j$ for $\omega\neq 0$, then
there exists a positive $t\in (0,1)$ such that
$F(\v)\le (F_i(\v-t\omega)+F_j(\v+t\omega))/2$.
\end{itemize}
If each $F_i$ is strictly convex, then $F$ is strictly convex.
\end{lemma}

Note that $\omega$ can be chosen arbitrary small, the condition in Lemma \ref{Piece} is 
a local property around $D_i\cap D_j$. 

\proof
The condition is clearly necessary. We prove the sufficiency.  Note that since $Z$ is of dimension $d-2$
or less, if the condition is valid for $\v\in (D_i\cap D_j)\setminus Z$ then it is
also valid for $D_i\cap D_j$ by continuity of convex functions.
Let us show the simplest case that $\R^d=D_1\cup D_2$
and $D_1\cap D_2$ is a hyperplane.  Take $\x\in D_1,\v\in D_2$ and find $\w\in [\x,\v]\cap D_1\cap D_2$.
By the assumption, if $0\neq \omega$ is a positive multiple of $\v-\x$, there exists $t>0$ \Red{such that}
\begin{equation}
\label{s}
F(\w)\le \frac 12 (F_1(\w-t \omega)+F_2(\w-t\omega)).
\end{equation}
and $\w-t \omega\in (\x,\w)$ and $\w+t \omega\in (\w,\v)$. Therefore we find $\mu_1,\mu_2\in (0,1)$
that $\w-t \omega=(1-\mu_1)\x+\mu_1 \w$ and $\w+t\omega=(1-\mu_2)\w+\mu_2 \v$. Using convexity of $F_i$, 
we have
\begin{equation}
\label{l}
F_1(\w-t\omega)\le (1-\mu_1) F_1(\x)+\mu_1 F(\w)
\end{equation}
and
\begin{equation}
\label{r}
F_2(\w+t\omega)\le (1-\mu_2) F(\w)+ \mu_2 F_2(\v).
\end{equation}
Using (\ref{s}),(\ref{l}),(\ref{r}) we deduce
$$
F(\w)\le \frac {1-\mu_1}{1-\mu_1+\mu_2} F_1(\x)+\frac {\mu_2}{1-\mu_1+\mu_2} F_2(\v).
$$
Because we can take arbitrary small $t$, 
the required convexity inequality holds for all $\x\in D_1,\v\in D_2$ and $\w\in (\x,\v)\cap D_1 \cap D_2$. 
Take $\u,\x\in D_1, \w\in D_1\cap D_2, \v\in D_2$ so that 
$\x,\w$ are within the open segment $(\u,\v)$. 
Take $\lambda,\mu\in (0,1)$ \Red{such that} $\w=(1-\lambda)\x+\lambda \v$ and $\x=(1-\mu)\u + \mu \w$. 
By the above discussion, we have
\begin{equation}
\label{12}
F(\w)\le (1-\lambda)F_1(\x)+\lambda F_2(\v).
\end{equation}
By the convexity of $F_1$,
\begin{equation}
\label{1}
F_1(\x)\le (1-\mu)F_1(\u)+\mu F(\w).
\end{equation}
From (\ref{12}) and (\ref{1}), we obtain
$$
F_1(\x)\le \frac {1-\mu}{1-\mu+\mu \lambda} F_1(\u)+\frac {\mu \lambda}{1-\mu+\mu\lambda} F_2(\v).
$$
Summing up, we know that any pair of two points $\u\in D_1$ and $\v\in D_2$,
the required convexity inequality is valid for any point $\x\in (\u,\v)$.
Therefore we can merge domains of convexity and 
the proof for the case $\R^d=D_1\cup D_2$ is finished. 
If each $F_i$ is strictly convex, then the resulting inequality is strict. 
One can easily extend this discussion to the general case, we simply 
repeat the merging process for adjacent domains sharing a \Red{codimension} one face. 
The set $Z$ does not disturb this merging process because $\{D_i\}$ 
are chain connected by
the adjacency relation induced by \Red{codimension} one faces. 
\qed
\medskip

A convex body $X$ is {\it strictly convex}, if $\x,\y\in X$ with $\x\neq \y$, then 
$(1-\lambda)\x+\lambda \y\in \mathrm{Inn}(X)$ for $\lambda\in (0,1)$, 
where $\mathrm{Inn}(X)$ is the interior of $X$. It is easy to see 
that a non empty set of the form $\{\v|\ F(\v)\le h\}$ for some $h>0$ 
is strictly convex if $F$ is strictly convex except $Y$ with 
a convex $Y\subset \mathrm{Inn}(X)$.

\begin{theorem}
\label{Strict}
$S$ is strictly convex.
\end{theorem}

\proof
As in the proof of Theorem \ref{Triangle}, considering $p_n$ as a variable $\v$, 
the surface $\partial(S)$ is a contour of the sum of $V_2(q_i,q_{i+1},\v)$
where $[q_i,q_{i+1}]\ (i=0,\dots, \ell-1)$ are the related visible edges. Here $q_i\in \{p_1,\dots, p_{n-1}\}$ and the index $i$ is considered modulo $\ell$. 
The hyperplanes which contain a face of $\Delta(p_1,\dots,p_{n-1})$ gives a partition $\{D_i\}_{i\ge 1}$ 
of $\R^d$ into a finite number of 
polyhedral cells and the set of visible faces is invariant within each $D_i$ outside 
$\Delta(p_1,\dots,p_{n-1})$. Let $F_i$ be the function
$A_3(p_1,\dots,p_{n-1},\v)$ restricted to $D_i$, and define a constant function
$F_0(\mathbf{v})=A_3(\Delta(p_1,\dots,p_{n-1}))$ for $\v\in D_0:=\Delta(p_1,\dots,p_{n-1})$.
Then $F_i(\v)=F_j(\v)$ for $\v\in D_i\cap D_j$ is clear. 
Let $Z$ be as in Lemma \ref{Piece} which is a finite set of $\R^3$.
We claim that the condition of Lemma \ref{Piece} is also satisfied.
Indeed in the same way as in the proof of Theorem \ref{Triangle},
$A_3(\Delta(p_1,\dots,p_{n-1},\v))$ is totally differentiable at $\v\in (D_i\cap D_j)\setminus Z$ 
with $1\le i<j$, i.e.,
\begin{equation}
\label{Tangent}
F(\u)=F(\v)+\nabla F(\v) \cdot (\u-\v)+o(\|\u-\v\|)
\end{equation}
with $\nabla F(\v) \neq (0,0,0)$.
Take $\omega\neq 0$ such that $\v-\omega\in D_i$, $\v+\omega\in D_j$.  If $F_i(\v-t_0 \omega)<F(\v)+\nabla F(\v)\cdot (-t_0 \omega)$ 
for some $t_0\in (0,1)$, then 
$(F_i(\v-t \omega)-F_i(\v))/t\le  (F_i(\v-t_0 \omega)-F_i(\v))/t_0\le M $ for all $t\in (0,t_0)$ with 
a constant $M<\nabla F(\v) \cdot (-\omega)$
by convexity of $F_i$. This contradicts (\ref{Tangent}) 
and we see $F_i(\v-t_0\omega)\ge F(\v)+\nabla F(\v)\cdot (-t_0\omega)$. 
In the same way, we have $F_j(\v+t_0 \omega)\ge F(\v)+\nabla F(\v)\cdot (t_0 \omega)$ and thus
$F(\v)\le (F_i(\v-t_0\omega)+F_j(\v+t_0\omega))/2$.
For $\v\in D_0\cap D_j$ with $j\ge 1$, the condition is trivial because $A_3(\Delta(p_1,\dots, p_{n-1},\v))\ge
A_3(\Delta(p_1,\dots, p_{n-1}))$. It remains to show that each $F_i\ (i\ge 1)$ 
is strictly convex to apply Lemma \ref{Piece}.

Clearly $\ell\ge 3$.
We claim that $V_2(q_i,q_{i+1},\x)$ is a convex function. Indeed, consider a plane $P_i$ passing $q_i$
perpendicular to $[q_i,q_{i+1}]$ and the orthogonal projection $g$ to $P_i$. Then we have
$V_2(q_i,q_{i+1},\x)=\| g(\x)-q_i\| \|q_{i+1}-q_i\|/2$. Since $g$ is linear and $g(q_i)=q_i$, 
\Blue{triangle} inequality implies
\begin{equation}
\label{Conv2}
\|g((1-\lambda) \x+\lambda \y)-q_i\|\le (1-\lambda)\|g(\x)-q_i\| + \lambda \|g(\y)-q_i\|
\end{equation}
for $\lambda\in [0,1]$ which proves the claim. As the sum of convex function is convex, we know
$\sum_{i=0}^{\ell-1} V_2(q_i,q_{i+1},\x)$ and $F$ are convex.
The equality for $\lambda\in (0,1)$ in (\ref{Conv2}) occurs only if $g(\x)-q_i$ and $g(\y)-q_i$ are
linearly dependent. 
This happens only when $\x,\y,q_i,q_{i+1}$ are in the same plane. 
However we can find an index \Red{$i$ such that} $\x,\y,q_i,q_{i+1}$ are not in the same plane.
Indeed, by our implicit assumption on visibility,  
the $\ell+2$ points $\{\x,\y\} \cup \{q_i|\  i=0,\dots, \ell-1\}$
\Red{cannot} be in the same plane.  Therefore, we always have
$$
\sum_{i=0}^{\ell-1} V_2(q_i,q_{i+1},(1-\lambda)\x+\lambda \y)
<(1-\lambda) \sum_{i=0}^{\ell-1} V_2(q_i,q_{i+1},\x)+\lambda \sum_{i=0}^{\ell-1} V_2(q_i,q_{i+1},\x)
$$
for $\lambda\in (0,1)$. This proves that $F$ is strictly convex except $\Delta(p_1,\dots,p_{n-1})$. 
Since $S=\{\v\in \R^3|\ F(\v)\le h\}$ for some $h>0$, we have shown the \Blue{theorem}.
\qed
\medskip

\begin{corollary}\label{Rigid}
A minimum $n$-hedron $\Delta_n=\Delta(p_1,\dots,p_n)$ does not allow deformation
of a single vertex, i.e., 
there exists a positive $r$ \Red{so that} if 
$\Delta(p_1,\dots,p_{n-1},\x)$ is a \Red{minimum} $n$-hedron with $\| \x-p_n\|<r$, 
then $\x=p_n$. 
\end{corollary}

\proof
Let $\Delta_n=\Delta(p_1,\dots,p_n)$ be the minimum $n$-hedron. By the proof of Theorem \ref{Triangle},
$p_n$ is on the boundary of the convex polyhedron $C$ as well as on the surface $\partial(S)$ which is a boundary
of the strictly convex set $S$ by Theorem \ref{Strict} and $S$ is contained in $C$. Take a small ball $B$ around
$p_n$ \Red{with the property that} $B\cap \partial(C)$ is contained in a single face of $C$.
If $B$ contains a point $\v\in \partial(C) \cap \partial(S)$ other than $p_n$, then the segment $[\v, p_n]$ is in $\partial(C)\cap S$ by convexity.
However since $S$ is strictly convex, $(\v+p_n)/2 \in \mathrm{Inn}(S)$ which contradicts $S\subset C$.
\qed

\section{\Red{Possible shapes of minimum $n$-hedron} for $n\le 12$}
\Red{By numerical experiments and Gr\"obner bases computation, one can 
give a list of possible shapes of minimum $n$-hedron for $n\le 12$.}

\begin{lemma}
\label{Incenter}
Let $X$ be a tetrahedron of vertices $K,L,M,N$
and $g$ be the orthogonal projection to the plane $P$ containing
 $L,M,N$.  
Let $K$ move in the plane parallel to $P$, keeping its
 volume $V_3(X)$ invariant. Among such $K$, the minimum surface 
 area $A_3(X)$
 is attained when $g(P)$ is the \Red{incenter} of the triangle $LMN$.
\end{lemma}

\proof
Let $h_1, h_2, h_3$ be the \Red{signed} height\footnote{\Red{It is positive in direction to the interior of the triangle $\Delta(L,M,N)$.}}
of the point $g(K)$ from the edge
$MN$, $NL$, $LM$ respectively in the plane $P$ and $h$ is the length of
the segment $[K,g(K)]$. 
Denote by $e_1,e_2,e_3$ the length of the edge 
$MN$, $NL$, $LM$ respectively. Then we have  
$V_2(L,M,N)=(e_1 h_1+ e_2 h_2+ e_3 h_3)/2$ and if
$V_2(L,M,N)$ is fixed, $(h_1,h_2,h_3)$ 
gives a coordinate system of points in $P$ under this constraint, 
i.e., two of $\{h_1,h_2,h_3\}$ determine the remainder through this relation.
Our problem is to minimize
$$
A_3(K,L,M,N)=
V_2(L,M,N)+ \frac12\left(e_1 \sqrt{h_1^2+h^2}+ e_2 \sqrt{h_2^2+h^2}+ e_3 \sqrt{h_3^2+h^2}
\right)
$$
under  $V_2(L,M,N)=(e_1 h_1+ e_2 h_2+ e_3 h_3)/2$. Since 
$|h_i|\rightarrow \infty$ for some $i$
implies $A_3(K,L,M,N)\rightarrow \infty$, we may assume
that $(h_1,h_2,h_3)$ are in a compact set of $\R^3$. Therefore
the minimum of $A_3(K,L,M,N)$ exists. Using Lagrange multiplier, we see that
the minimum is attained when
$$
\frac{\partial}{\partial h_i}
\left(
A_3(K,L,M,N)-\lambda (V_2(L,M,N)-(e_1 h_1+ e_2 h_2+ e_3 h_3)/2) \right) =0,
$$
for $i=1,2,3$. This implies $h_1/\sqrt{h_1^2+h^2}=h_2/\sqrt{h_2^2+h^2}=h_3/\sqrt{h_3^2+h^2}$ and consequently $h_1=h_2=h_3$. 
Therefore the minimum is attained when $g(P)$ 
is the \Red{incenter} of the triangle $LMN$.
\qed

\Red{This may be a well-known result. 
One of the referees informed me of a similar discussion in an encyclopedia on elementary
geometry, "Kikagaku Dai Jiten" vol. 5, page 440, ed. Shikou Iwata, in Japanese.}

\medskip
\Red{
Similarly if a pyramid $X$ whose base $k$-gon $B$ is circumscribed about a circle and 
its apex $K$ moves in the plane parallel to $B$, 
then $A_3(X)$ is minimized when the orthogonal projection of the apex to $B$ is the center of the circle. In fact,
let us define $e_i, h_i$ in a similar manner. Though lengths $h_i$ are determined by two parameters, say $h_1$ and $h_2$, we minimize the surface area function
in a less constrained domain $$\left\{ 0\le h_i\le C\left|\ \sum e_i h_i=2 A_3(B)\right.\right\}\subset \R^k$$
with a sufficiently large $C>0$. Then the condition $h_1=h_2=\dots=h_k$ 
is attained at the center of the circle
under the assumption.
}

\begin{lemma}
\label{Radon}
Any $d+2$ points in $\R^d$ is partitioned into two non empty disjoint sets $U$ and $V$ \Red{such that} $\Delta(U)\cap \Delta(V)\neq \emptyset$.
\end{lemma}

\proof  This is due to Radon (\cite[Theorem 1.1.5]{Schneider:93} or 
\cite{Gruenbaum:03}). It is an easy consequence of the
linear dependence of $v_i-v_{d+2}$ for $i=1,\dots, d+1$ for any point set
$\{v_1,v_2\dots, v_{d+2}\}$.
\qed
\medskip

A {\it $k$ bi-pyramid} is a polygon composed of two 
pyramids sharing the same $k$-gon base joined base to base. 
A {\it regular $k$ bi-pyramid} is a bi-pyramid composed of two congruent regular pyramids
sharing the regular $k$-gon base. 
Its main diagonal is the segment joining two apexes vertically
passing the center of the base.

\begin{lemma}
\label{Bi}
Among $k$ bi-pyramids $\Delta$ whose convex bases are circumscribed about a circle of radius $h$,
the minimum 
$$
\frac{A_3(\Delta)}{V_3(\Delta)^{2/3}}=\left(3^{7/2}k\cot\left(\frac {(k-2)\pi}{2k}\right)\right)^{1/3} 
$$
is attained when it is a regular bi-pyramid whose main diagonal has length \Red{$\sqrt{8} h$}.
\end{lemma}

\proof
The minimization
problem is divided into two pyramids, say, an upper pyramid and a lower pyramid. 
Let $B$ be the common base polygon. Letting $\theta_i \ (i=1,\dots,k)$ be the \Red{vertex} angles of $B$, we obtain 
$A_2(B)=hr$ with $r= 2 \sum_{i=1}^k \cot(\theta_i/2)$ and $V_2(B)=h A_2(B)/2$.
Let $H_1$ be the height of the apex of the upper pyramid $\Delta'$
to the base $B$, and $H_2$ is the one for the lower pyramid $\Delta''$.
Then we have $V_3(\Delta')=V_2(B) H_1/3= h^2 r H_1/6$. \Red{By the discussion after Lemma \ref{Incenter}, 
$A_3(\Delta')$ is minimized when the orthogonal projection of the apex to $B$ is the center of the circle
and}
$$
A_3(\Delta')-V_2(B)=\frac 12 A_2(B) \sqrt{h^2+H_1^2}=\frac 12 hr \sqrt{h^2+H_1^2}.
$$
Let us fix $r$ and minimize \Red{$A_3(\Delta)$ by selecting $H_1$, $H_2$ and $h$
keeping $V_3(\Delta)$} invariant. 
\Red{Fixing $h^2 (H_1+H_2)/6$, the minimum of $h(\sqrt{h^2+H_1^2}+\sqrt{h^2+H_2^2})/2$ 
is attained when $H_1=H_2=\sqrt{2}h$.} 
Now we have
$A_3(\Delta)= \sqrt{3} h^2r$ and $V_3(\Delta)=\Red{\sqrt{2}h^3 r /3}$. Thus
$$
\frac {A_3(\Delta)}{V_3(\Delta)^{2/3}} = 3^{7/6} 2^{-1/3} r^{1/3}.
$$
Since $\cot(x/2)$ is convex for $x\in (0,\pi)$, by Jensen's inequality, 
the minimum of $r$ is achieved by the regular $k$-gon when $\theta_1=\theta_2=\dots=
\theta_k=\pi-2\pi/k$ and $r=2k \cot(\frac {(k-2)\pi}{2k})$.
\qed
\medskip

\begin{theorem}
\label{Minimum}
We have\footnote{For the case $n=6$, see \S \ref{Prob} (2).}
\begin{eqnarray*}
\alpha_4&=& 6\cdot 3^{1/6}\approx 7.20562,\\
\alpha_5&=&3^{5/3}\approx 6.24025 ,\\
\alpha_6&\le& \eta_6=3^{7/6}2^{2/3}\approx 5.71911,\\
\alpha_7&\le& \eta_7=3^{7/6} 5^{5/12} (\sqrt{5}-2)^{1/6}\approx 5.53841,\\ 
\alpha_8&\le& \eta_8\approx 5.42118, \\
\alpha_9&\le& \eta_9\approx 5.31637,\\
\alpha_{10}&\le& \eta_{10}\approx 5.2533\\
\alpha_{11}&\le& 5.20713\\
\alpha_{12}&\le &  \eta_{12}=3^{7/6} (70 - 30 \sqrt{5})^{1/3}\approx 5.14835
\end{eqnarray*}
where $\eta_8, \eta_9,\eta_{10}$ are algebraic numbers of degree $72$, $78$, $36$ 
respectively. $\alpha_4$ is attained by a regular tetrahedron 
and $\alpha_5$ by a regular $3$ bi-pyramid.
\end{theorem}

Our experiments suggest that all the inequalities are equalities,
though we did not identify the exact value for $\alpha_{11}$. 
Several specialists working on computer science told me that  
brute force optimization does not seem feasible as it has too much free variables for now.

\proof
\Red{The case $n=4$ may belong to a folklore. At least a written proof is found 
in Hadwiger \cite[p.273, (187)]{Hadwiger} using \Blue{the} Steiner symmetrization.  Here we give a direct proof.}
Let $KLMN$ be the \Red{minimum $4$-hedron}. 
By Lemma \ref{Incenter},
projection of $K,L,M,N$ to the corresponding basis triangle must be its \Red{incenter}.
Let $KH$ be the perpendicular from $K$ to $\Delta(L,M,N)$ and $HI$, $HJ$ be the perpendicular
from $H$ to $LM$ and $LN$. Since $H$ is the \Red{incenter} of $\Delta(L,M,N)$, we have $HI=HJ$, 
$KI\bot LM$, $KJ\bot LN$.
From $KI=KJ$ and $IL=JL$, we see that $\angle KLM=\angle KLN$. By cyclic
discussion we see, $\angle KLM=\angle KLN=\angle MLN:=\angle L$. Similarly we see, three angles at each vertex of
$\Delta(K,L,M,N)$ are identical for all vertices, which are denoted by $\angle K, \angle L, \angle M, \angle N$.
Since the sum of angles of triangular faces are all equal to $2\pi$, we deduce
that $\angle K=\angle L=\angle M=\angle N$, therefore all the faces are regular triangles. This proves the case of the minimum $4$-hedron. 

For \Red{minimum} $5$-hedron, 
in light of  Lemma \ref{Decrease} we may assume that none of vertices 
is contained in the convex hull of remaining four vertices. 
Therefore by Lemma \ref{Radon},  five vertices are divided into two sets $\{K,L,M\}$ and $\{N,O\}$ \Red{for which} $\Delta(K,L,M)\cap \Delta(N,O)\neq \emptyset$.  
The problem is therefore reduced to Lemma \ref{Bi} for $k=3$.
This case was also shown in \cite[Theorem 5.5]{Boroczky2_Iso:97}.

For $n\ge 6$, we performed a random search of the minimum. 
A rough sketch of the empirical method is 

\begin{enumerate}
\item Choose random $n$ points in $\R^3$ and determine the combinatorial structure of the convex hull, in particular the {\it valency vector}, that is, 
the multi-set of valencies of vertices. 
\item Iterate process 1, until we find a valency vector of small variance. 
Experimentally, 
we know that $A_3(\Delta)/V_3(\Delta)^{2/3}$ \Red{cannot} be small if this variance is large.
\item Select a vertex, an edge or a face of \Red{$\Delta$ and minimize $A_3(\Delta)/
V_3(\Delta)^{2/3}$} by moving its extremities, keeping the valency vector invariant. 
If the valency vector changes, then we skip this minimization.
\item Find two points $v_1,v_2$ which gives the diameter of $\Delta$ 
and apply an affine transformation to make a little smaller the diameter  but 
keeping the plane orthogonal to $v_1-v_2$ invariant. 
\item Repeat several times these processes 2,3 and 4 at random.
\end{enumerate}

Until $n\le 12$ it seems the above iteration leads us to a possible minimum for a fixed valency vector. Trying many valency vectors, we can guess the target shape. 
Then we perform algebraic computation
to obtain the exact \Blue{minimum} configuration. Taking into account the expected symmetry 
of the target shape, we set up a system of 
algebraic equations with a small number of variables. Then we eliminate 
variables by using some program equipped with Gr\"obner basis computation. We used 
Mathematica, PARI-GP and Risa-Asir appealing to each advantage. Gr\"obner basis 
computation has a lot of subtleties. Successful computation depends heavily
on the number of variables, their imposed order, and degree of polynomials. 
Hereafter we describe our computation but skipping such technical details, giving 
necessary information to reconfirm the computation.

\begin{figure}[htbp]
\begin{center}
\subfigure{
\includegraphics[height=6cm]{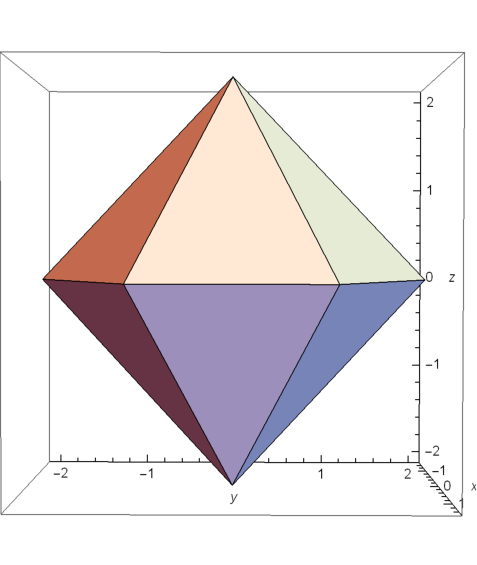}
}
\qquad
\subfigure{
\includegraphics[height=6cm]{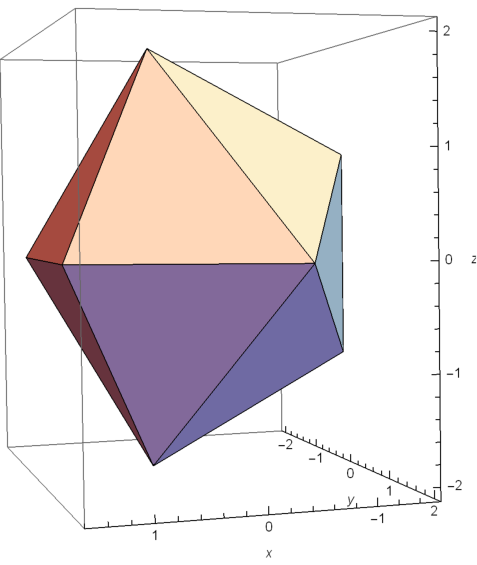}
}
\qquad
\subfigure{
\includegraphics[height=6cm]{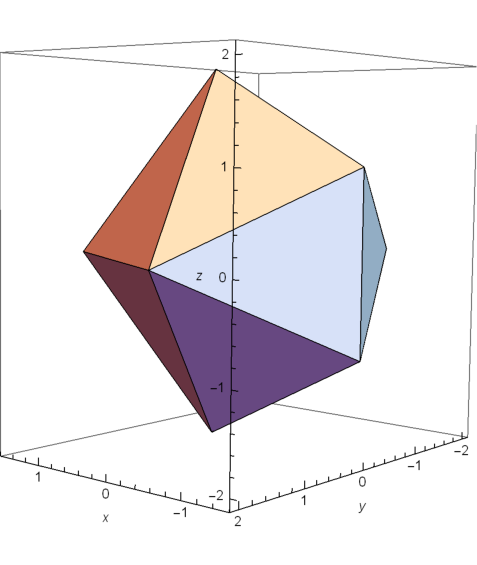}
}
\qquad
\subfigure{
\includegraphics[height=6cm]{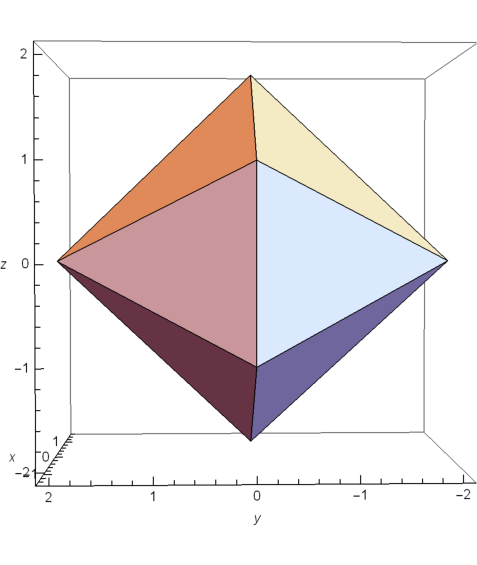}
}
\end{center}
\caption{Minimum $8$-hedron}
\label{Shape}
\end{figure}

\begin{figure}[htbp]
\begin{center}
\includegraphics[height=20cm]{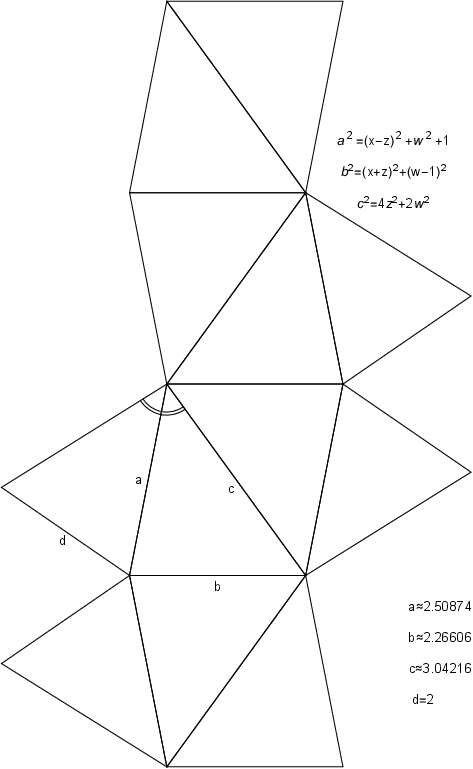}
\end{center}
\caption{Origami diagram}
\label{Dev}
\end{figure}

By our experiments, 
the target shapes for $\eta_6$ and $\eta_7$ are attained by regular bi-pyramid as in Lemma \ref{Bi}.
The most difficult and interesting shape appears when $n=8$, see \Blue{Figures} \ref{Shape} and \ref{Dev}. 
\Red{It is combinatorially equivalent to the Siamese dodecahedron, one of the deltahedra.}

\Red{By using numerical minimization, we} could guess that
$8$ points are of the form:
$$
(z,0,\pm w), (x,\pm 1,0), (-z,\pm w, 0), (-x,0,\pm 1)
$$
with
$$
w\approx 2.0428, x \approx 1.53525, z \approx 0.476614.
$$
\Red{Let us consider $w,x,z$ as variables and obtain their exact algebraic representations.}
We have
$$
A_3(\Delta)=4\sqrt{w^2 + (x - z)^2} + 4\sqrt{(w-1)^2 w^2 + w^2 (x + z)^2 + (w(x - z) + 2z)^2}
$$
and 
$$
V_3(\Delta)=4w(x + wx + z)/3. 
$$
We view $A_3/V_3^{2/3}$ as the function of three variables. Basically
our task is to eliminate valuables 
from 
$$
\partial_x \left(\frac{A_3^3}{V_3^2}\right)=\partial_w \left(\frac{A_3^3}{V_3^2}\right)=\partial_z \left(\frac{A_3^3}{V_3^2}\right)=0.
$$
First transfer problems into the one on polynomials with integer coefficients, \Red{by} putting
$$
u=\sqrt{w^2 + (x - z)^2}, v=\sqrt{(w-1)^2 w^2 + w^2 (x + z)^2 + (w(x - z) + 2z)^2}.
$$
Then eliminate $u,v$ to find an ideal over $x,w,z$
and perform primary ideal decomposition (this was indispensable for this computation).
We obtain \Red{minimal} polynomials\footnote{\Red{Minimal polynomials of $x$ and $z$ are the ones subject to substitution
$t\rightarrow t^2$.}} of $w,x^2,z^2$:
\begin{align*}
&8 - 40t - 32t^2 + 268t^3 - 14t^4 + 378t^5 - 916t^6 + 874t^7\\
& - 265t^8 - 314t^9 + 374t^{10} - 150t^{11} + 21t^{12},
\end{align*}
\begin{eqnarray*}
&&-1500625 + 246891400t - 6498924184t^2 + 197676252320t^3 - 549916476544t^4 \\
&&+ 9593743607488t^5 - 37068998078592t^6 +  43451585720832t^7 + 6412940883200t^8 \\
&&- 47369088623616t^9 + 34505601388544t^{10} - 10887830962176t^{11} + 1413638553600t^{12},
\end{eqnarray*}
\begin{eqnarray*}
&&-881721 + 14088624t - 507815656t^2 + 22228266304t^3 - 345876361600t^4 \\
&& + 2163078191936t^5 - 5229062814592t^6 +  2885777661952t^7 + 604100406528t^8 \\
&& + 284044459008t^9 - 1111813844992t^{10}+ 65086242816t^{11} + 157070950400t^{12}.
\end{eqnarray*}
The \Red{minimal} polynomial of $A_3(\Delta)^6/(V_3(\Delta))^{4}$ is
\begin{eqnarray*}
&&846253032058341803633618097683156083357246027504784634537836544 \\
&-&  145765911302088136407360046924472940590350227969907327078760448t \\
&+&  44739094836549297939345827315732094525400511681413644681216t^2 \\
&-& 5444218664651134627342263572192722894788633799480098816  t^3 \\
&+& 381929202246269536064619254896305729053865712762224t^4 \\
&-& 23215968331655851588483378342178431615039134384 t^5 \\
&+&  908544689594387775769635417411363042641304 t^6 \\
&-& 26376155703404842068063899980163109720 t^7 \\
&+&  639590587552165626186327476412759 t^8 \\
&-& 9114814042610279966292752064 t^9 + 144758783681628174471168 t^{10} \\
&-&  130494391161126912 t^{11} + 4980736000 t^{12}
\end{eqnarray*}
A non-trivial coincidence of two angles indicated in Figure \ref{Dev} is \Red{confirmed by algebraic computation
of cosine values of the angles}.
One can also confirm numerically that this minimum shape is rigid, see Section \ref{Prob}.
\medskip

\begin{figure}[htbp]
\begin{center}
\subfigure[$n=9$]
{
\includegraphics[height=5cm]{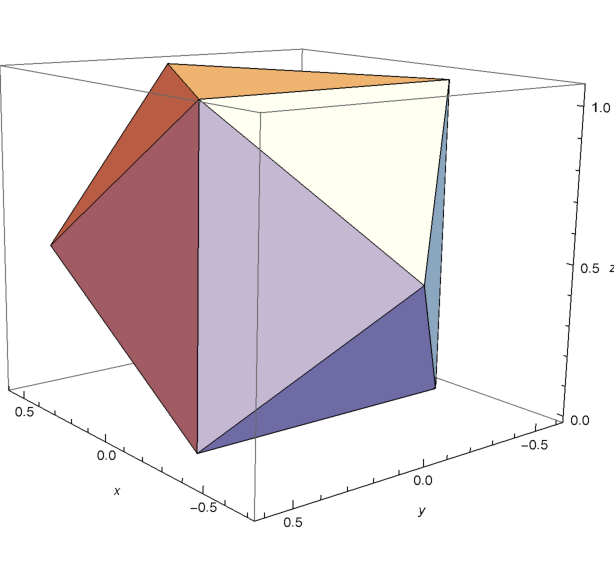}
}
\qquad
\subfigure[$n=10$]
{
\includegraphics[height=5cm]{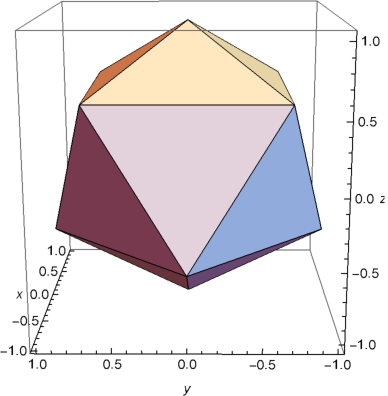}
}
\end{center}
\caption{Minimum polyhedron\label{9_10}}
\end{figure}

For $\eta_9$, consider a regular triangular prism 
and put three identical $4$-pyramids to
each of rectangular side faces 
whose centroid is the foot of the perpendicular from the apex of the pyramid, see Figure \ref{9_10} (a).
Let the edge length of the regular triangle be $1$. Then the height of the prism $b$, and
the height of the $4$-pyramid $h$ are expected to be
$$
b\approx 1.04725, h \approx 0.413823.
$$
We have
$$
A_3(\Delta)=\frac{\sqrt{3}}2 + 3b \sqrt{h^2 + \frac14} + 3 \sqrt{h^2 + \frac{b^2}4}
$$
and
$$
V_3(\Delta)=\frac{\sqrt{3}b}4 + b h.
$$
We treat $A_3(\Delta)/V_3(\Delta)^{2/3}$ as a function of two variables $b$ and $h$
and apply the elimination 
of variables as we did in $n=8$. Note that to treat $\sqrt{3}$, we also 
introduce another variable $s$ and the polynomial $s^2-3$ to be added in the ideal. 
The \Red{minimal} \Blue{polynomials} of $b^2$ and $h^2$ are
\begin{eqnarray*}
&-&3600-9384 t+157415 t^2+1871849 t^3-3005515 t^4-3048555 t^5+7100157 t^6\\
&-&716904 t^7-5370867t^8+3887865 t^9-810945 t^{10}-53622 t^{11}+17415 t^{12}+2187 t^{13},
\end{eqnarray*}
and 
\begin{eqnarray*}
&-&27-216 t-5688 t^2+99268 t^3+2629424 t^4-11859776 t^5-198587904 t^6\\
&+&641098752 t^7+2269974528 t^8+3790651392 t^9-43985534976 t^{10}\\
&+&74140876800 t^{11}-37371248640 t^{12}+5435817984 t^{13}.
\end{eqnarray*}
The \Red{minimal} polynomial of $A_3^6/V_3^4$ is
\begin{eqnarray*}
&-&8741200275671730192755167246352564248392781977833773782269952 \\
&+&5692272790315788765597663433429575175625193067065671949484032t \\
&+& 1663401637275489431763071207791450034909825698981382756499456t^2 \\
&+& 205636897183575223972130099822721877708248944269405343514624t^3 \\
&-& 29496333327693613396843751515776856015704029599701614592t^4 \\
&-& 58714195329202332973530206007453465620049796957569024t^5 \\
&-& 13596161545396297014562622838466932898374596846592t^6 \\
&-& 847586880386300377059351613641377507384112384t^7 \\
&-& 58377287904203791631778906263194550638656t^8 \\
&+& 3993703760487214498878732921512576256t^9\\
&-& 12610065164386918027558684269276t^{10} \\
&+& 4281392126518694452576397473t^{11} - 20704119330241635606528t^{12}\\
&+&  21761395104153600t^{13}
\end{eqnarray*}

For $\eta_{10}$, prepare an anti-prism,  a convex hull of a square and its parallel
square rotated by $\pi/4$, and put two identical regular $4$-pyramids on the two parallel squares, see 
Figure \ref{9_10} (b).
We introduce a coordinate of $10$ points: 
$$
(\pm 1,0,-h),(0,\pm 1,-h),(\pm \frac1{\sqrt {2}},\pm \frac1{\sqrt{2}},h),(0,0,\pm z)
$$
with
$$
h \approx 0.541397, z \approx 1.02619
$$
and minimize
$$
\frac{A_3^3}{V_3^2}
=\frac{36 \left(\sqrt{3 - 2\sqrt{2} + 8h^2} + 
    \sqrt{1 + 2h^2 - 4hz + 2z^2}\right)^3}{(h + \sqrt{2}h + z)^2}.
$$
The \Red{minimal} polynomials of $h^2,z^2,A_3^6/V_3^4$ are
$$
1+48 t+144 t^2-16128 t^3-31296 t^4+273408 t^5+28672 t^6,
$$
$$
47089+1130960 t-1729392 t^2+2846464 t^3-1889856 t^4-277504 t^5+28672 t^6,
$$
and
\begin{eqnarray*}
&-& 9592639401335565227088041861971968 +  362253880325110957404812476416t \\
&-&  4924615865029090098020352t^2 +  462296427139672731648t^3 \\
&-& 713296009601244t^4 +  274678452t^{5} + t^{6}.
\end{eqnarray*}

\begin{figure}[htbp]
\begin{center}
\includegraphics[height=5cm]{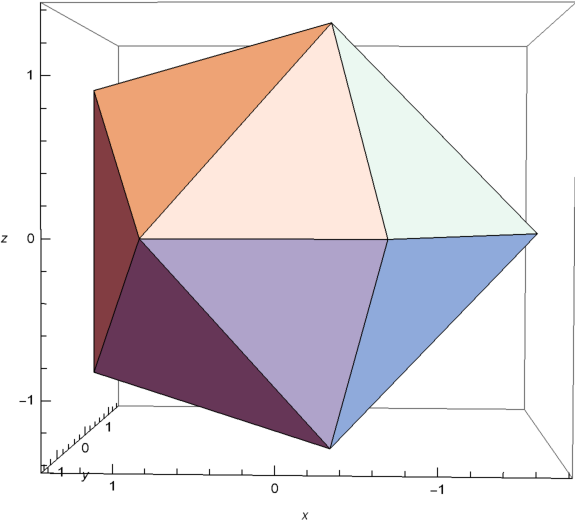}
\end{center}
\caption{$n=11$: From $y$-axis direction}
\label{11P}
\end{figure}

\begin{figure}[htbp]
\begin{center}
\subfigure[From $x$-axis \Blue{positive} direction]{
\includegraphics[height=4.5cm]{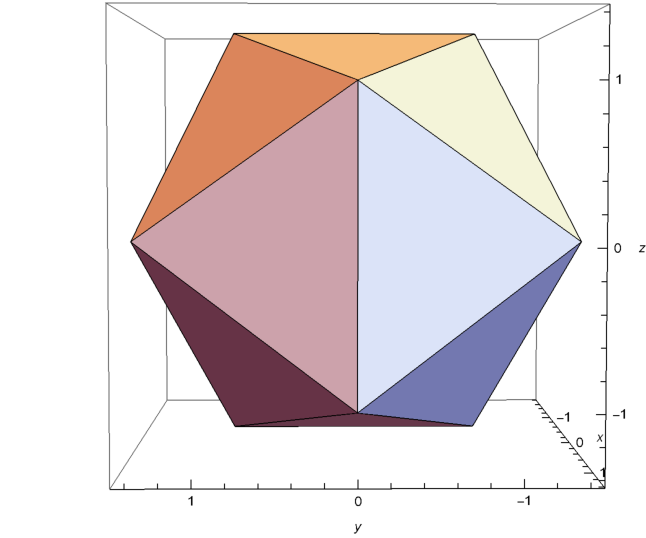}
}
\qquad
\subfigure[From $x$-axis \Blue{negative} direction]{
\includegraphics[height=4.5cm]{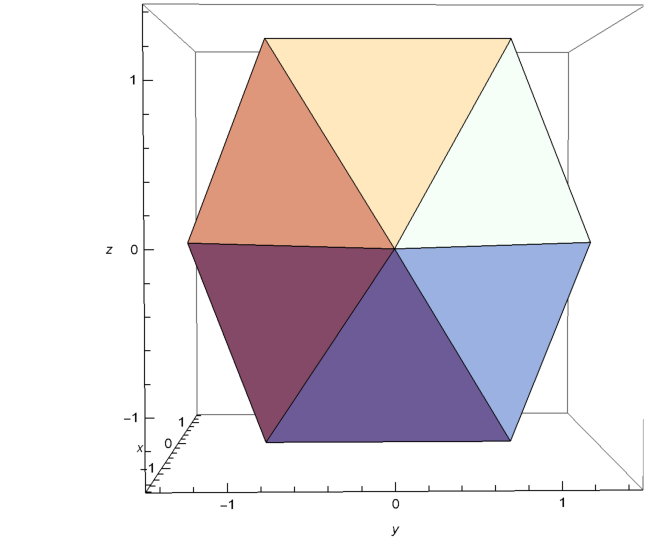}
}
\end{center}
\caption{Minimum 11-hedron\label{11y}}
\end{figure}

We also obtained the conjectural shape for $n=11$ by experiments. It is a convex hull of 
$$
(x_1, 0, \pm 1),
(x_2, \pm y_1,0), 
(-x_3, \pm y_2,0), 
(-x_4, \pm y_3,\pm z),(-x_5, 0, 0)
$$
with
$$
x_1 \approx 1.15135, x_2 \approx 0.617047, x_3 \approx 0.91681, x_4 \approx 0.550702, 
  x_5 \approx 1.98113,
$$
$$
y_1 \approx 1.4264, y_2 \approx 1.34059,  y_3 \approx 0.845054, z \approx 1.38959
$$
with
$$
A_3/V_3^{2/3}\approx 5.207134373504469,
$$
see \Blue{Figures} \ref{11P} and \ref{11y}. \Red{It has 18 faces, which follows from Euler's formula and 
the fact that all faces are triangles.} \Blue{The shape is combinatorially equivalent to a polyhedron obtained by merging two adjacent vertices
of the regular icosahedron into one.}

We could not make the \Red{coordinates} algebraic, because the expected symmetry group $(\Z/2\Z)^2$ 
is too small, and the number of valuables is too large.

The minimum $12$-hedron is of course expected to be the regular icosahedron with
$$
A_3/V_3^{2/3}= 3^{7/6} (70 - 30 \sqrt{5})^{1/3}\approx 5.14835.
$$
\section{Problems}\label{Prob}

We give a list of intriguing problems. 

\begin{enumerate}
\item Can we give an asymptotic estimate for the convergence of $(\alpha_n)$ ?
\item Prove our candidates \Blue{minimum} for $i=6,7,8,9,10,12$. (Added in Revision: 
the validity for $i=6$ is confirmed in \cite{Boroczky_Kovacs:19}).
\item Is minimum $n$-hedron $\Delta(p_1,\dots,p_n)$ rigid ?  
We say that 
$\Delta(p_1,\dots,p_n)$ is rigid if it does not allow deformation of $n-3$ vertices, i.e., 
there exists a positive $r$ \Red{so that} for any subset $I$ of $\{1,\dots,n\}$ of 
cardinality $n-3$, if $\Delta(x_1,\dots,x_n)$ is a minimum $n$-hedron
with $\|x_i-p_i\|<r$ for $i\in I$ and $x_i=p_i$ for $i\not \in I$, then $x_i=p_i$ holds for all $i$.
\item Is the symmetry group of the 
minimum $n$-hedron non-trivial for all $n$~? Can it have a chirality, i.e. , can its symmetry group in $O(3)$
and that in $SO(3)$ be different ?
\item Is $\partial(S)$ in Theorem \ref{Strict} defined by strongly convex function ?
\end{enumerate}
\bigskip

{\bf Acknowledgement.}
The author wishes to express his deepest gratitude to anonymous referees 
for pointing out mistakes, missing references and very careful reading of the original manuscript.


\end{document}